# Weak units and homotopy 3-types

André Joyal and Joachim Kock

*For Ross Street, on his 60th birthday*

Abstract. We show that every braided monoidal category arises as End($I$) for a weak unit $I$ in an otherwise completely strict monoidal 2-category. This implies a version of Simpson's weak-unit conjecture in dimension 3, namely that one-object 3-groupoids that are strict in all respects, except that the object has only weak identity arrows, can model all connected, simply connected homotopy 3-types. The proof has a clear intuitive content and relies on a geometrical argument with string diagrams and configuration spaces.

## 0. Introduction

The subtleties and challenges of higher category theory start with the observation (in fact, not a trivial result) that not every weak 3-category is equivalent to a strict 3-category. The topological counterpart of this is that not every homotopy 3-type can be realised by a strict 3-groupoid. The discrepancy between the strict and weak worlds can be pinpointed down to the case of connected, simply-connected 3-types, where it can be observed rather explicitly: such 3-types correspond to braided monoidal categories (in fact braided categorical groups), while connected, simply-connected strict 3-categories are essentially commutative monoidal categories — the braiding is forced to collapse, as a consequence of the Eckmann-Hilton argument. In precise terms, strict $n$-groupoids can realise only homotopy $n$-types with trivial Whitehead brackets.

This collapse can be circumvented by weakening the structures. The notion of tricategory of Gordon, Power, and Street [**1**] is meant to be the weakest possible definition of 3-category. They show that a tricategory with only one object is equivalent to a Gray monoid, and in particular, a tricategory with one object and one arrow is equivalent to a braided monoidal category. Furthermore, every braided monoidal category arises in this way. The most general result relating higher categories to homotopy types is Tamsamani's theorem [**9**], that weak $n$-groupoids (in the sense of Tamsamani) can realise all homotopy $n$-types. This result was conjectured by Grothendieck [**2**], or rather: it was stated as a desideratum for any future theory of weak higher categories.

In Tamsamani's theory, and in most other theories of higher categories, the essential weakening bears on the composition laws and their interchange laws. However, a careful analysis of the situation in strict 3-groupoids led Simpson [**8**] to observe that the collapse of the braiding, via the Eckmann-Hilton argument, can be traced back to the strictness of the identity arrows. He conjectured that (a suitable notion of) strict $n$-groupoids with weak identity arrows should realise all homotopy $n$-types, and furthermore that the homotopy category of such $n$-groupoids should be equivalent to the homotopy category of Tamsamani $n$-groupoids. (In fact he went further and conjectured that the same homotopy equivalence should hold in the non-invertible case, i.e. for $n$-categories, not just for $n$-groupoids.) An ad hoc notion of weak identity arrows was sketched, but the details were not worked out, and it was acknowledged that it might not be the correct notion to fulfil the conjectures. Simpson's conjectures are highly surprising: they go against all trends in higher category theory, where the emphasis is mostly on the composition laws, and questions about identity arrows are often swept under the carpet. A consequence of the conjectures is that every weak $n$-category should be equivalent to one with strict composition laws and strict interchange laws!

In this work we prove a version of Simpson's conjecture in the crucial case of dimension 3. We restrict ourselves to the connected, simply-connected case, working with strict monoidal 2-categories with weak units. The basics of weak units in monoidal 2-categories is worked out in a companion paper [**3**], but in fact very little is needed in our proof. Our key result is this:

**Main Theorem.** *Let $I$ be a weak unit of an otherwise completely strict monoidal 2-category. Then $\operatorname{End}(I)$ is a braided monoidal category, and every braided monoidal category arises in this way.*

Connected, simply-connected homotopy 3-types correspond to braided categorical groups. Under the correspondence of the Main Theorem, these correspond to strict 2-groupoids with invertible tensor product and weak units, which in turn can be regarded as one-object 3-groupoids. Hence we get the following version of Simpson's conjecture in dimension 3:

**Main Corollary.** *One-object 3-groupoids that are strict in all respects, except that the object has only weak identity arrows, can model all connected, simply connected homotopy 3-types.*

The paper is organised as follows. In Section 1 we show that $\operatorname{End}(I)$ is braided, and explain the geometry of this braiding. In Section 2 we introduce the geometric language of train track diagrams and show that the space of all train track diagrams is acyclic. Finally in Section 3, given a braided monoidal category $\mathbb{B}$, we use a geometrical construction to get a monoidal 2-category with weak unit $I$ such that $\operatorname{End}(I)$ is equivalent to $\mathbb{B}$.

## 1. From weak unit to braiding

**1.1. Semi-monoidal 2-categories.** A strict semimonoidal 2-category (or a 2-category with strict multiplication) is a (strict) 2-category $\mathscr{C}$ equipped with a strictly associative multiplication functor $\otimes : \mathscr{C} \times \mathscr{C} \to \mathscr{C}$. We write the tensor product by plain juxtaposition: $(X, Y) \mapsto XY$. We use the symbol # to denote composition of arrows, written from the left to the right, writing for example $f \# g$ for the composite

$$\xrightarrow{f} \xrightarrow{g}$$

and we use the same symbol for 'horizontal' composition of 2-cells. We use the symbol Ⓒ to denote identity 2-cells.

**1.2. Weak units,** cf. [**3**]. A *weak unit* in $\mathscr{C}$ is a pair $(I, \alpha)$ where $I$ is an object in $\mathscr{C}$ with the property that tensoring with $I$ from either side is an equivalence of 2-categories $\mathscr{C} \to \mathscr{C}$, and $\alpha : II \to I$ is an equi-arrow in $\mathscr{C}$ (i.e., an arrow admitting a quasi-inverse $I \to II$).

It is shown in [**3**] that this notion of weak unit is equivalent to the definition that can be extracted from the notion of tricategory of Gordon, Power, and Street [**1**] involving usual left and right constraints $IX \to X \leftarrow XI$, and it is also equivalent to the notion of weak unit that can be extracted from the abstract notion of fair categories of [**6**], which is a general 'non-algebraic' approach where the emphasis is on the contractible space of all units, not on any arbitrary fixed unit itself.

The key point for these results, and all we need to know for the present purposes, is that there is a canonical 2-cell $\mathsf{A} : I\alpha \Leftrightarrow \alpha I$. For the reader's convenience we briefly outline the construction, referring to [**3**] for all details. Since tensoring with $I$ on the left is an equivalence of 2-categories, for each object $X$ the functor $\mathrm{Hom}(IX, X) \to \mathrm{Hom}(IIX, IX)$ is an equivalence of categories, and in particular, essentially surjective. Hence the essential inverse image of $\alpha X$ is non-empty, so we can choose an arrow $\lambda_X : IX \to X$ together with an invertible 2-cell $I\lambda_X \Leftrightarrow \alpha X$. Similarly, tensoring instead with $I$ on the right, we can find an arrow $\rho_X : XI \to X$ and an invertible 2-cell $X\alpha \Leftrightarrow \rho_X I$. Now take $X = I$, and paste together four instances of the chosen 2-cells as follows:

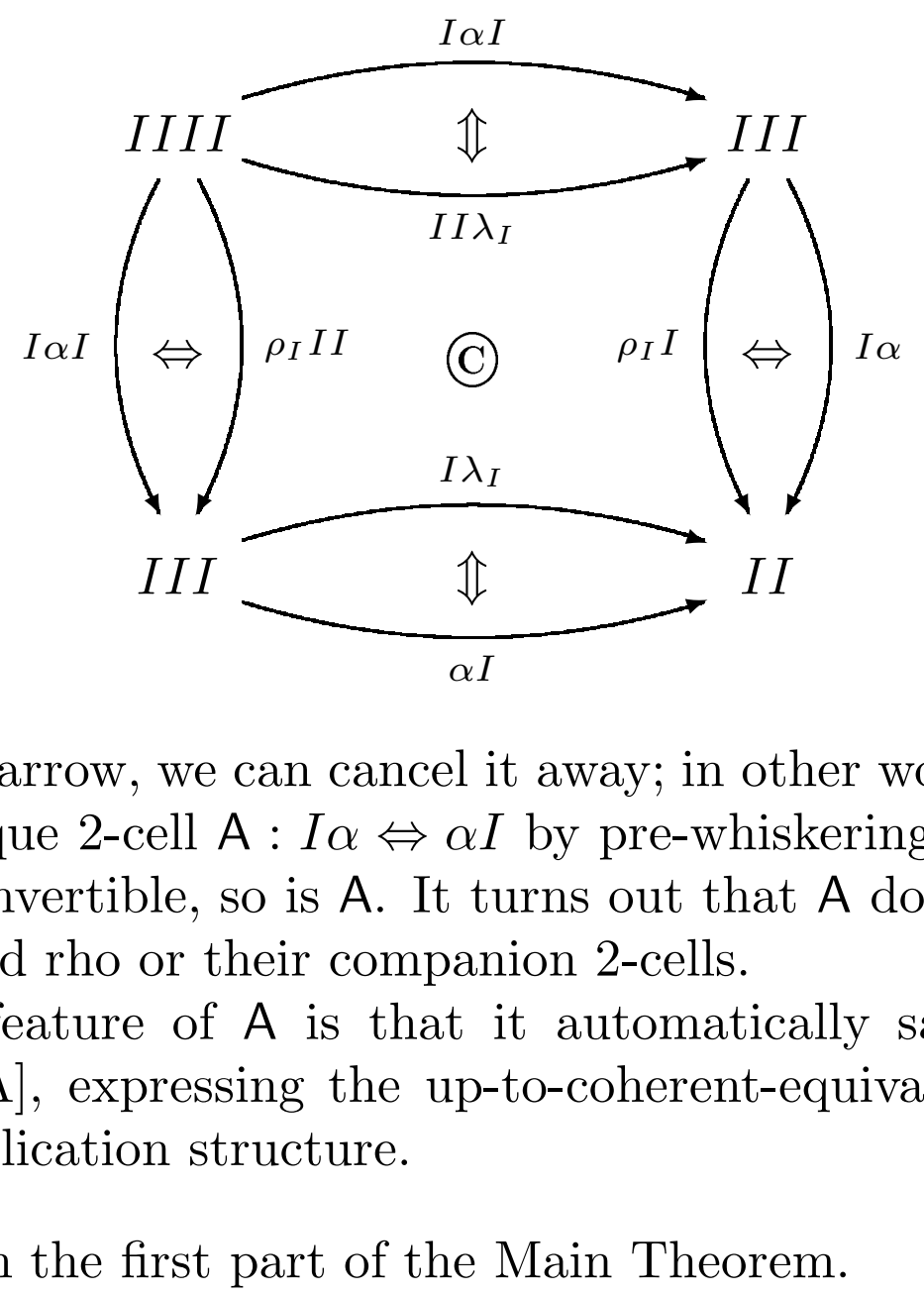

Since $I\alpha I$ is an equi-arrow, we can cancel it away; in other words the total 2-cell is obtained from a unique 2-cell $\mathsf{A} : I\alpha \Leftrightarrow \alpha I$ by pre-whiskering with $I\alpha I$. Since the involved 2-cells are invertible, so is $\mathsf{A}$. It turns out that $\mathsf{A}$ does not depend on the choices of lambda and rho or their companion 2-cells.

The surprising feature of $\mathsf{A}$ is that it automatically satisfies the pentagon equation [**3**, Thm. A], expressing the up-to-coherent-equivalence associativity of $\alpha$ viewed as a multiplication structure.

We now establish the first part of the Main Theorem.

Proposition 1.3. *Let* $(I, \alpha)$ *be a weak unit in* $\mathscr{C}$*. Then the strict monoidal category* $(\mathrm{End}(I), \#, \mathrm{id}_I)$ *is naturally braided.*

The situation resembles that of a strict monoid in ***Gray***: given $f$ and $g$ in $\mathrm{End}(I)$ we need to provide an invertible 2-cell

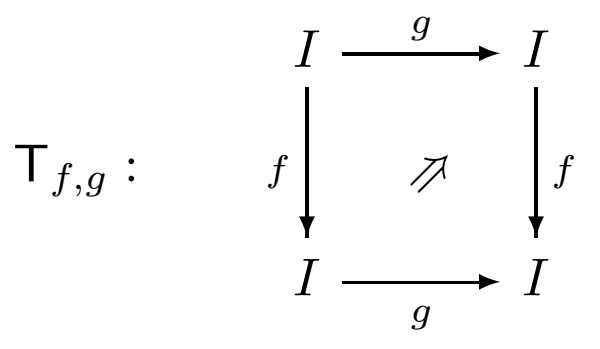

and check the axioms for a braiding.

Lemma 1.4. *There is an invertible* 2*-cell in* $\mathscr{C}$*:*

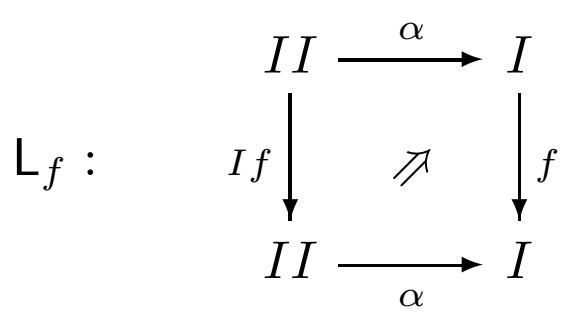

*natural in* $f$*.*

Proof. $\mathsf{L}_f$ is defined as the unique 2-cell satisfying the equation

(1) 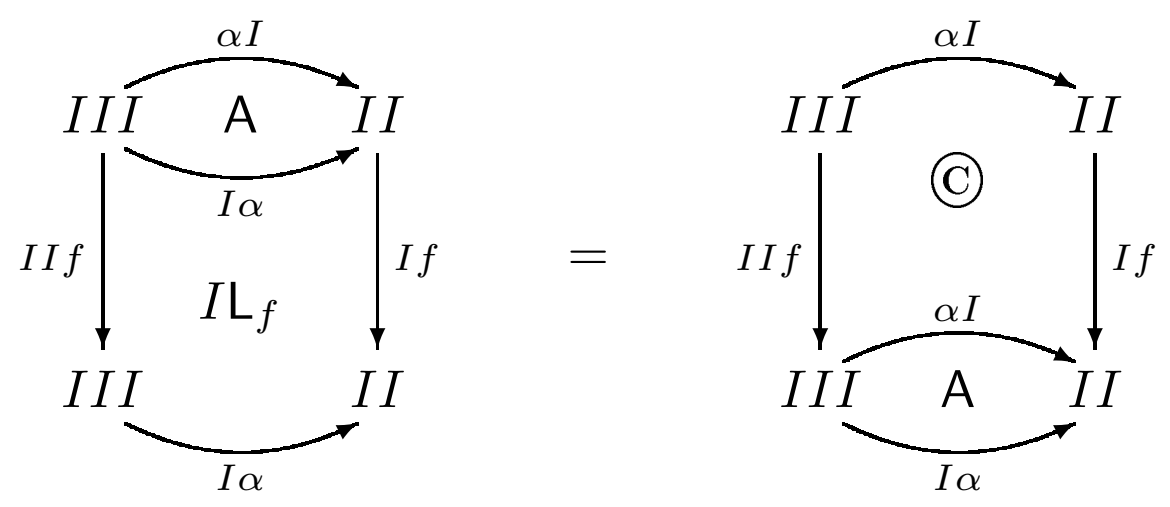

This makes sense: since the three other 2-cells in the diagram are invertible, the cell labelled $I\mathsf{L}_f$ is well-defined, and since tensoring with $I$ on the left is an equivalence of 2-categories and hence a bijection on the level of 2-cells, also $\mathsf{L}_f$ itself is uniquely defined.

Naturality in $f$ means that for any 2-cell $\mathsf{U} : f \Rightarrow g$, we have

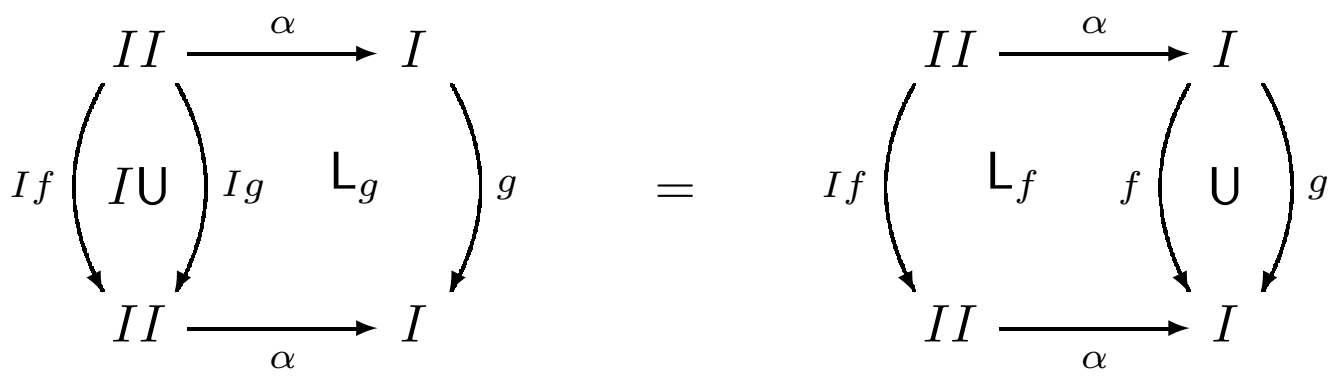

To check this equation, it is enough to check it holds after tensoring with $I$ on the left. This allows us to use the defining property (1) of $\mathsf{L}$ on both sides, and then the result follows from the trivial observation that this equation holds:

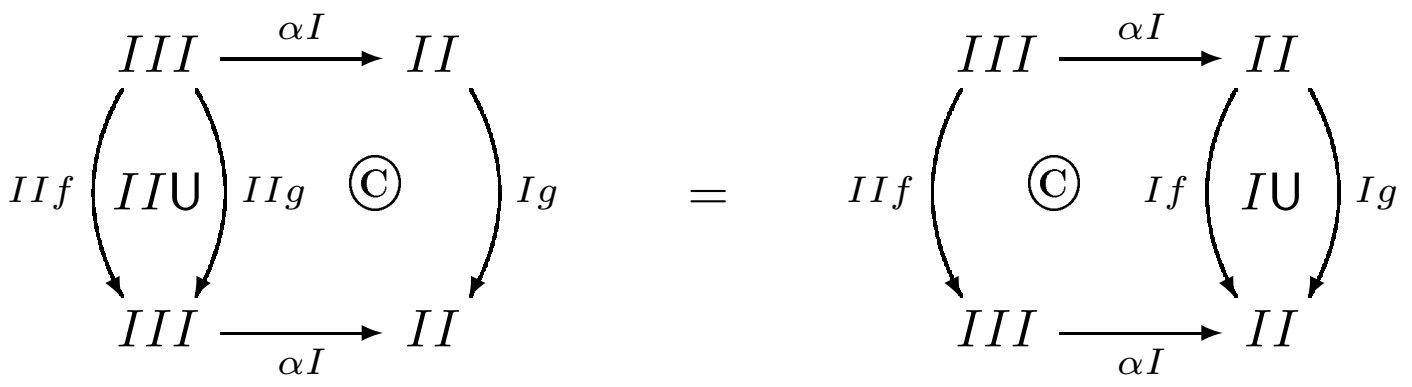

□

Remark 1.5. There is another description of $\mathsf{L}_f$, given in [**3**]: modulo a canonical isomorphism $\lambda_I \Leftrightarrow \alpha$ it is an instance of a naturality 2-cell for the left constraint whose construction was outlined in 1.2.

Lemma 1.6. *The* 2*-cells* $\mathsf{L}_f$ *are compatible with composition of endomorphisms of* $I$*:*

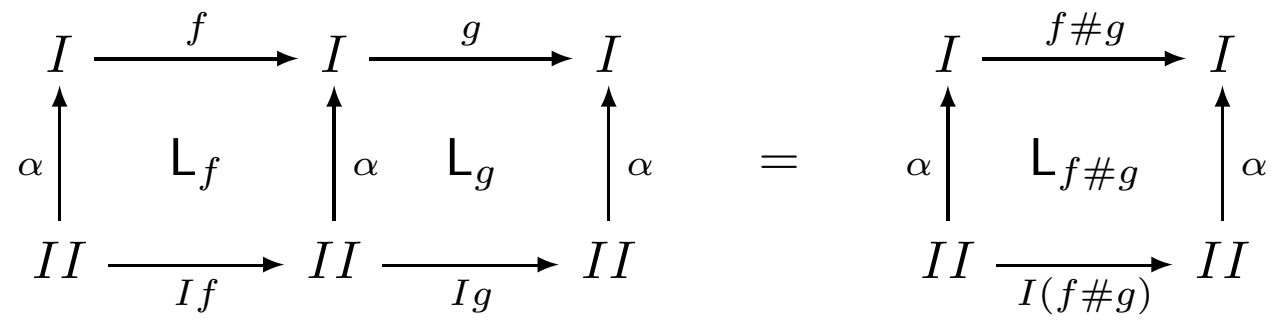

*Also, if* $f$ *is the identity arrow of* $I$*, then* $\mathsf{L}_f$ *is the identity* 2*-cell.*

Proof. After tensoring with $I$ on the left, the left-hand side of the equation is computed by gluing the two defining cylinder diagrams together along their common base $\mathsf{A}$. The result is clearly the defining cylinder for $\mathsf{L}_{f\#g}$. The statement about identity arrows also follows immediately from the defining cylinder diagram. □

**1.7. Mates.** Let $\beta : I \to II$ be a right adjoint to $\alpha$, with counit $\mathsf{E} : \beta \# \alpha \Rightarrow \mathrm{id}_I$. We will abusively draw this 2-cell as

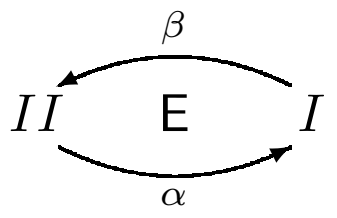

confident that the reader will remember that the source is $\beta \# \alpha$ and the target $\mathrm{id}_I$. There is a natural 2-cell $\mathsf{B} : \beta I \Leftrightarrow I\beta$ defined by the following equation of 2-cells from $\beta I \# I\alpha$ to $\mathrm{id}_{II}$:

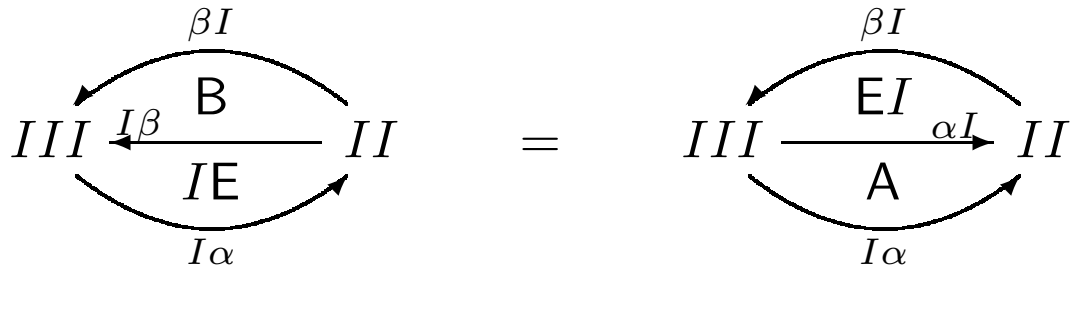

Using $\beta$ and $\mathsf{B}$ instead of $\alpha$ and $\mathsf{A}$, we get natural invertible 2-cells

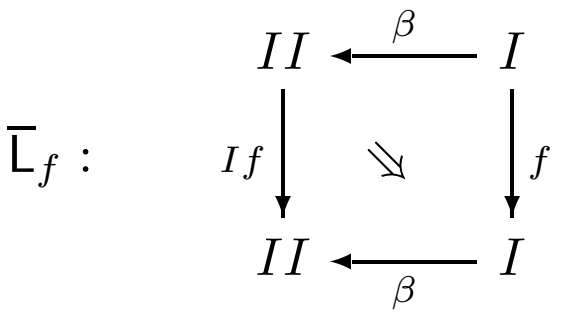

This is the mate of $\mathsf{L}_f$, cf. the following lemma. (Usually [**5**], mates are defined in terms of both the unit and counit, but it is practical for what follows to express the mate relation in terms of the counit only. This is possible since $\alpha$ is an equi-arrow.)

LEMMA 1.8. *We have this equation of* 2*-cells from* $\beta \# If \# \alpha$ *to* $f$*:*

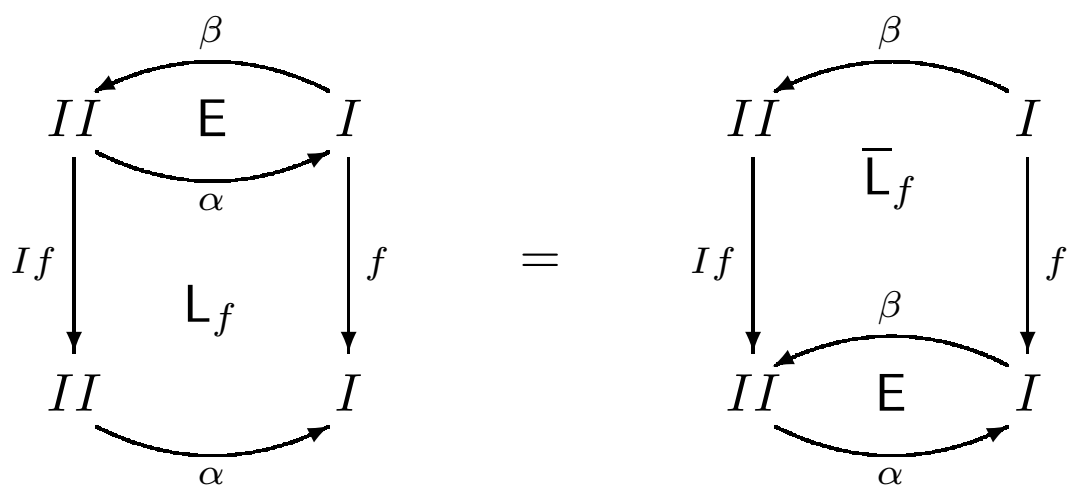

PROOF. Tensor with $I$ on the left and use the definition of $\mathsf{L}_f$, $\overline{\mathsf{L}}_f$, and $\mathsf{B}$. □

Finally, we will also need the corresponding constructions where instead the $I$-factor is on the right: there are natural invertible 2-cells

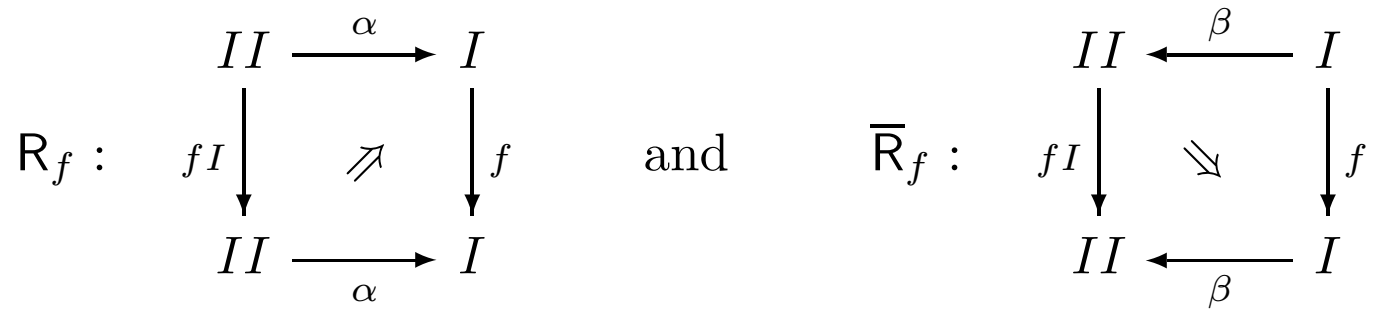

satisfying obvious analogues of Lemma 1.6 and Lemma 1.8.

PROOF OF PROPOSITION 1.3. The wanted 2-cell $\mathsf{T}_{f,g}$ is given by this pasting diagram

(2)

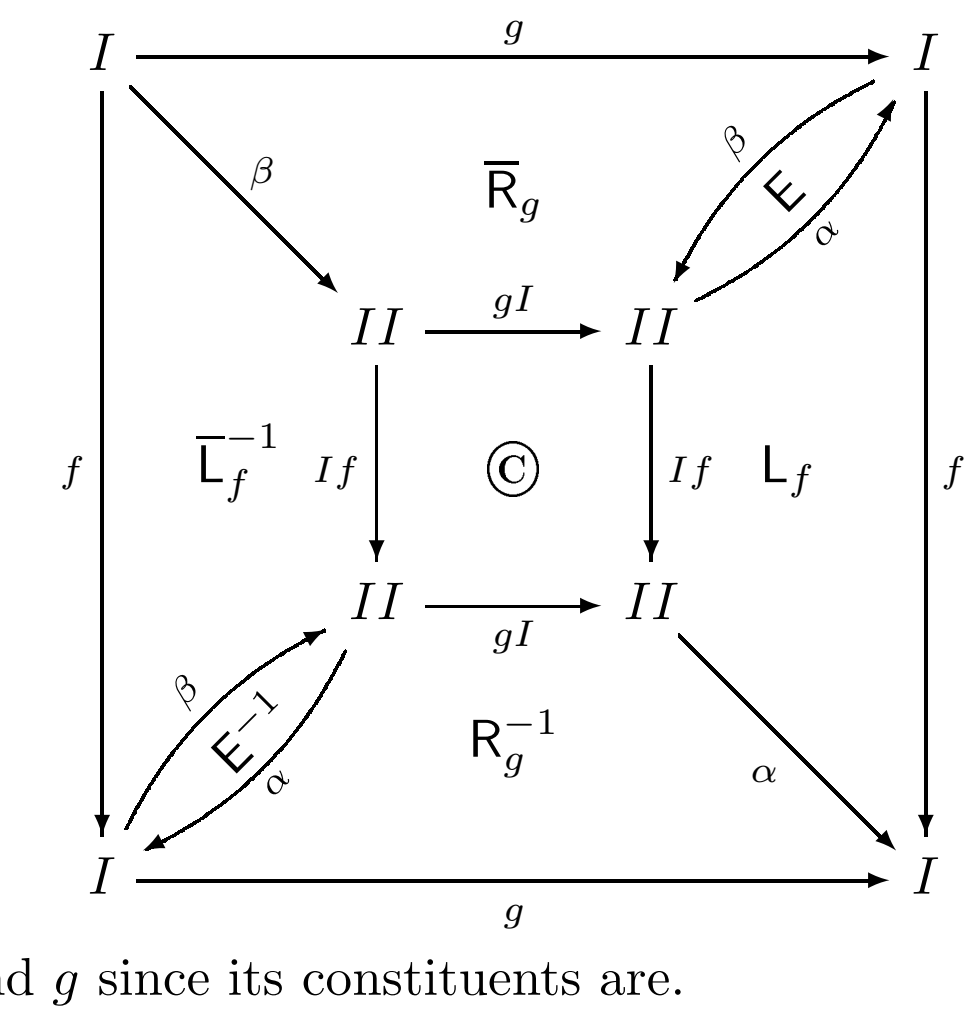

It is natural in $f$ and $g$ since its constituents are.

To see that the 2-cells $\mathsf{T}_{f,g}$ form a braiding we must check the triangle axioms, i.e., commutativity of the two diagrams

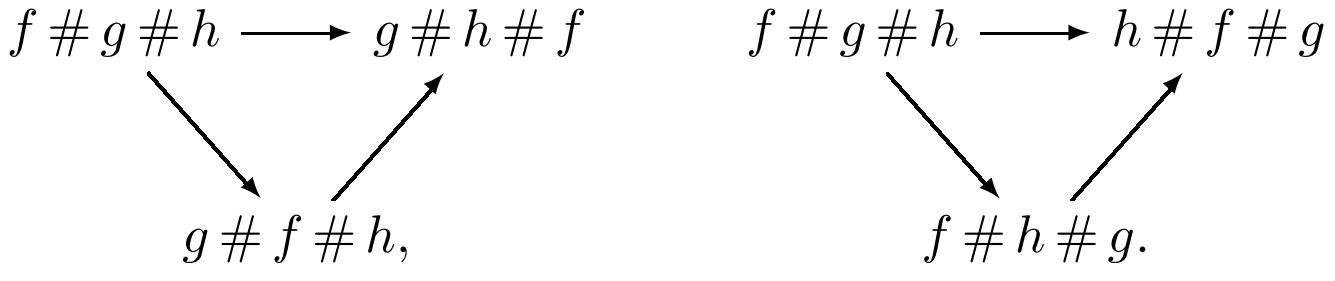

For the left-hand diagram, this means

$$\begin{array}{ccccc} I & \xrightarrow{\;g\;} & I & \xrightarrow{\;h\;} & I \\ {\scriptstyle f}\Big\downarrow & \mathsf{T}_{f,g} & \Big\downarrow{\scriptstyle f} & \mathsf{T}_{f,h} & \Big\downarrow{\scriptstyle f} \\ I & \xrightarrow[\;g\;]{} & I & \xrightarrow[\;h\;]{} & I \end{array} \quad = \quad \begin{array}{ccc} I & \xrightarrow{\;g\,\#\,h\;} & I \\ {\scriptstyle f}\Big\downarrow & \mathsf{T}_{f,g\#h} & \Big\downarrow{\scriptstyle f} \\ I & \xrightarrow[\;g\,\#\,h\;]{} & I. \end{array}$$

To establish this, spell out the diagrams in terms of (2), use Lemma 1.8 to cancel four cells near the middle of the diagram, and apply Lemma 1.6 twice. (The right-hand triangle axiom is checked using right-hand versions of 1.6 and 1.8.) □

REMARK 1.9. Note that $\mathsf{T}_{f,g}$ does not depend on the choice of $\beta$ used in its construction. This follows from essential uniqueness of adjoints: if an alternative $\beta'$ were used in the construction, the unique comparison 2-cell $\beta \Leftrightarrow \beta'$ would appear on the sides of $\overline{\mathsf{L}}_f^{-1}$ and $\overline{\mathsf{R}}_g$ and the net result would be the same.

**1.10. Geometry of the braiding.** In order to get an understanding of the nature of the braiding $\mathsf{T}_{f,g}$, a graphical interpretation is helpful. In fact, our proof that this braiding is generic will consist in taking these drawings literally.

Ignoring the weak unit structure of $\mathscr{C}$, it is a strict semi-monoidal 2-category, and as such it has an underlying strict semi-monoidal category whose algebraic structure can be expressed geometrically in terms of progressive planar string diagrams, in the usual way (cf. [**4**]; see also 2.1 below). The basic arrows appearing in the constructions are represented like this (to be read from the bottom to the top):

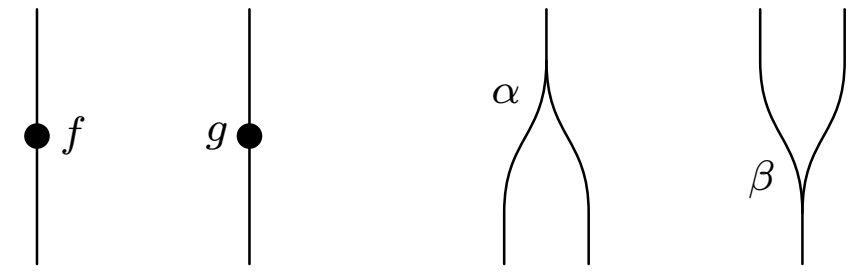

Each string represents a copy of $I$. An arrow is represented by a dot with some input strings coming in from the bottom (its source), and some output strings coming out at the top (the target). For reasons that will become clear, we suppress the dots for $\alpha$ and $\beta$. The tensor product of two arrows is represented by drawing their string diagrams side by side (parallel connection). Note that since the tensor product is not assumed to have a unit, the empty diagram is not permitted. Composition of arrows is realised by connecting the input strings of the second arrow to the output string of the first (serial connection).

The 2-cells do not have a proper geometric representation, but may be understood in terms of transformations of diagrams. For example, the 2-cell $\mathsf{L}_f$ is pictured like this:

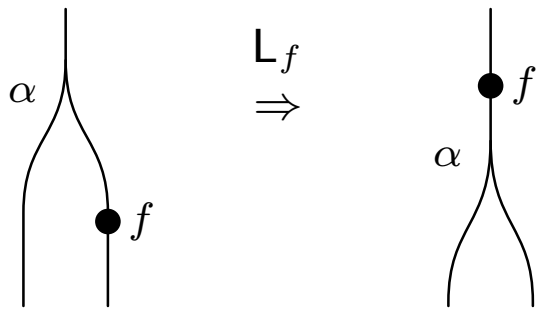

Now we can draw the sequence of seven 2-cells that make up $\mathsf{T}_{f,g}$:

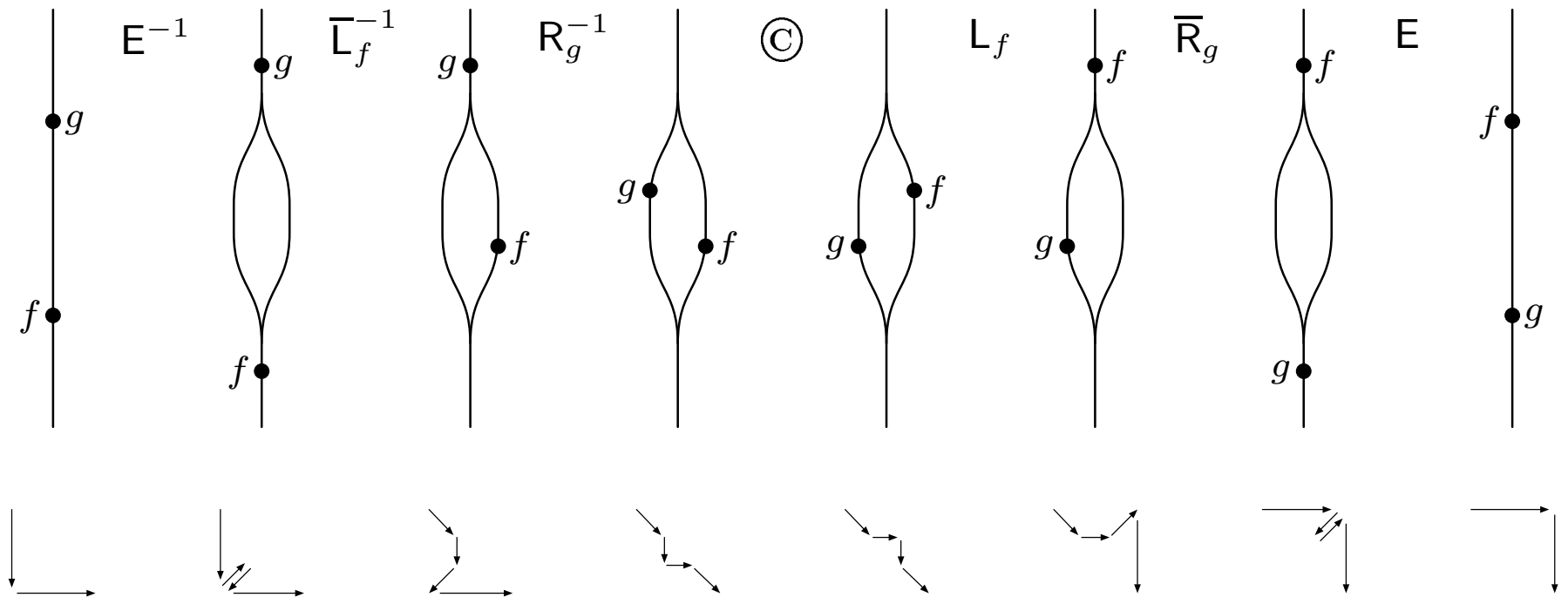

The small diagrams below the string diagrams refer to the corresponding path in the pasting diagram (2).

The point is that $f$ and $g$ change place and pass around each other in an orderly planar way, like two trains. The picture suggests that going left past each other is not the same as going right past each other, and that the braiding is not in general a symmetry. The result in Section 3 shows that indeed every braiding occurs in this way, and the proof consists in taking these diagrams seriously.

## 2. Train track diagrams

**2.1. Progressive plane diagrams.** We shall briefly recall some notions and results from *The geometry of tensor calculus, I* [**4**]. A *progressive plane graph* (between levels $b_0$ and $b_1$) is a finite graph $\Gamma$ (with boundary) explicitly embedded in $\mathbb{R} \times [b_0, b_1]$ such that

(i) the boundary of the graph is it intersection with $\mathbb{R} \times \{b_0, b_1\}$, and
(ii) the projection $\mathbb{R} \times [b_0, b_1] \to [b_0, b_1]$ is injective on each edge.

The vertices on level $b_0$ (resp. $b_1$) are called *inputs* (resp. *outputs*) of the graph; the remaining vertices are called *nodes*. Condition (ii) induces an orientation on each edge, and for each node an obvious notion of input and output edges of that node; the set of input edges and the set of output edges are both naturally ordered. A node has *valence* $(p, q)$ if it has $p$ input edges and $q$ output edges.

A *deformation* (or isotopy) of progressive plane graphs as above is a continuous function

$$h : \Gamma \times [0, 1] \to \mathbb{R} \times [b_0, b_1]$$

such that, for all $t \in [0, 1]$, the function

$$h(-, t) : \Gamma \to \mathbb{R} \times [b_0, b_1]$$

is a progressive plane graph (between levels $b_0$ and $b_1$).

A *progressive plane diagram* in a monoidal category $\mathscr{C}$ is a progressive plane graph whose nodes are labelled by arrows in $\mathscr{C}$, and whose edges are labelled by objects in $\mathscr{C}$, subject to the following compatibility condition: if a node $N$ is labelled by an arrow $f$, then the tensor product of the labels of the input edges of $N$ is the source of the arrow $f$, and the tensor product of the labels of the output edges of $N$ is the target of $f$.

A *deformation* of progressive plane diagrams in $\mathscr{C}$ is a deformation of the underlying progressive plane graphs whose labelling is constant on $\Gamma$.

The progressive plane diagrams (or just graphs) form a topological space (the topology is the compact-open topology for each fixed $\Gamma$). The deformations are the paths in this space.

The importance of progressive plane diagrams is that to each progressive plane diagram there is associated an arrow in $\mathscr{C}$, and this association is invariant under deformations (cf. [**4**], Thm. 1.2). Hence it makes sense to interpret drawings in $\mathscr{C}$ like in the previous section.

Data collections appropriate for generating monoidal categories are called tensor schemes in [**4**]. A tensor scheme $\mathscr{D}$ is the data of $\mathscr{D}_0^* \leftleftarrows \mathscr{D}_1$, where $\mathscr{D}_0$ is a collection of *objects*, $\mathscr{D}_0^*$ denotes the collection of all finite words in $\mathscr{D}_0$, and $\mathscr{D}_1$ is a collection of *arrows*, each having a source word and a target word.

Progressive diagrams make sense also in tensor schemes, and it is shown in [**4**, Thm. 1.3] that the free monoidal category on a tensor scheme $\mathscr{D}$ is the category whose objects are $\mathscr{D}_0^*$ and whose arrows are isotopy classes of progressive plane diagrams in $\mathscr{D}$.

**2.2. Effective diagrams and semi-monoidal categories.** If a progressive plane diagram in a monoidal category $\mathscr{C}$ has no inputs (resp. outputs), it must be interpreted as an arrow in $\mathscr{C}$ whose source (resp. target) is the unit object. In order to adapt the theory of [**4**] to semi-monoidal categories, clearly it is necessary to exclude nodes with empty in- or output: By an *effective* (progressive) plane graph we understand a non-empty progressive plane graph such that every node has at least one input edge and at least one output edge. Now the notion of an effective plane diagram in a semi-monoidal category $\mathscr{C}$ is obvious, and a semi version of [**4**, Thm. 1.2] follows.

**2.3. Train track diagrams.** A *train track diagram* in a semi-monoidal category is a progressive plane diagram such that every node has valence $(1,1)$, $(2,1)$, or $(1,2)$; in other words, an effective progressive plane diagram where no node has total valence greater than 3.

Effective diagrams, and in particular train track diagrams, are much more rigid than general progressive diagrams:

PROPOSITION 2.4. *The space $X_{\mathscr{C}}$ of all effective diagrams in $\mathscr{C}$ is acyclic, i.e. homotopy equivalent to $\pi_0(X_{\mathscr{C}})$. In other words, the space of diagrams deformation equivalent to a given effective diagram is contractible.*

PROOF. Ultimately, the reason is that each connected component of the complement of an effective diagram is contractible, which in turn is true because there are no 'islands': every part of the diagram is attached to the input and output levels. Here are the details:

It is enough to prove the proposition for the space of graphs $X$, since clearly $X_{\mathscr{C}} \to X$ is a covering projection, so we can forget about the labels. The proof is by induction on the number of edges. Let $X_\Gamma$ denote the space of graphs deformation equivalent to a given effective plane graph $\Gamma \subset \mathbb{R}^2$ (between levels $b_0$ and $b_1$). We shall construct a finite sequence of continuous maps

$$X_\Gamma = X_{\Gamma_0} \to X_{\Gamma_1} \to \cdots \to X_{\Gamma_n} = X_\emptyset$$

such that each map has contractible fibres. The last space is the singleton space consisting of the empty graph (between $b_0$ and $b_1$); each of the other spaces $X_{\Gamma_i}$ is a space of effective plane graphs deformation equivalent to some $\Gamma_i$. Each map will consist in either erasing a $(1,1)$-valent node or erasing an edge.

Erasing a $(1,1)$-node. Suppose $\Gamma$ contains a $(1,1)$-valent node $f$. Let $\Gamma'$ denote the shape obtained by erasing $f$: it is understood that the two edges adjacent to $f$ are joined to form a single edge $e$. Clearly $\Gamma'$ is effective if $\Gamma$ is. This operation of removing $f$ extends uniquely to a continuous map $\pi_f : X_\Gamma \to X_{\Gamma'}$. For any fixed graph $D'$ in $X_{\Gamma'}$, the $\pi_f$-fibre over $D'$ consists in all the possible ways of putting an extra node such that the resulting graph is deformation equivalent to $\Gamma$. These possibilities are parametrised by the inner points of the segment $e$, which is a contractible space.

Deleting an edge. We shall identify certain edges that can always be removed without spoiling the effectivity of the diagram. A *complete track* in a progressive plane graph $\Gamma$ is a sub progressive graph homeomorphic to a closed interval, going from an input of $\Gamma$ to an output of $\Gamma$. (So a non-empty progressive plane graph is effective when through every node there is a complete track.) An effective graph contains a *rightmost complete track*: it is the unique complete track with the property that there are no nodes or edges to the right of it. (Specifically, start with the rightmost input and progress: at each node you come to, turn right, choosing the rightmost output edge.) For general progressive plane graphs the notion of rightmost complete track is not well-defined, since there may be isolated connected components of the graphs floating around out to the right.

A *removable right edge* is an edge in the rightmost complete track such that if removed, the remaining graph is still a valid effective graph or possibly the empty graph. In other words, the edge does not start in an $(n,1)$-node and does not end in a $(1,n)$-node. (In the picture below, the only removable right edge is $(g_0, g_1)$.) Clearly the notions of rightmost complete track and removable right edge are invariant under deformation.

For a given graph $\Gamma$ with a chosen removable right edge $r$, let $\Gamma'$ denote the graph resulting from erasing that edge. The projection $X_\Gamma \to X_{\Gamma'}$ has contractible fibres: indeed, for a given graph $D'$ in $X_{\Gamma'}$ the possible ways of drawing a right edge from some node $g_0$ (at level $d_0$) to another node $g_1$ (at level $d_1$) is parametrised by the space of continuous functions on the interval $[d_0, d_1]$ dominating the rightmost track of $\Gamma'$, and with appropriate boundary values (to have $g_0$ and $g_1$ as endpoints). This space is clearly contractible. The area for the graph of such function is indicated in grey in the following figure.

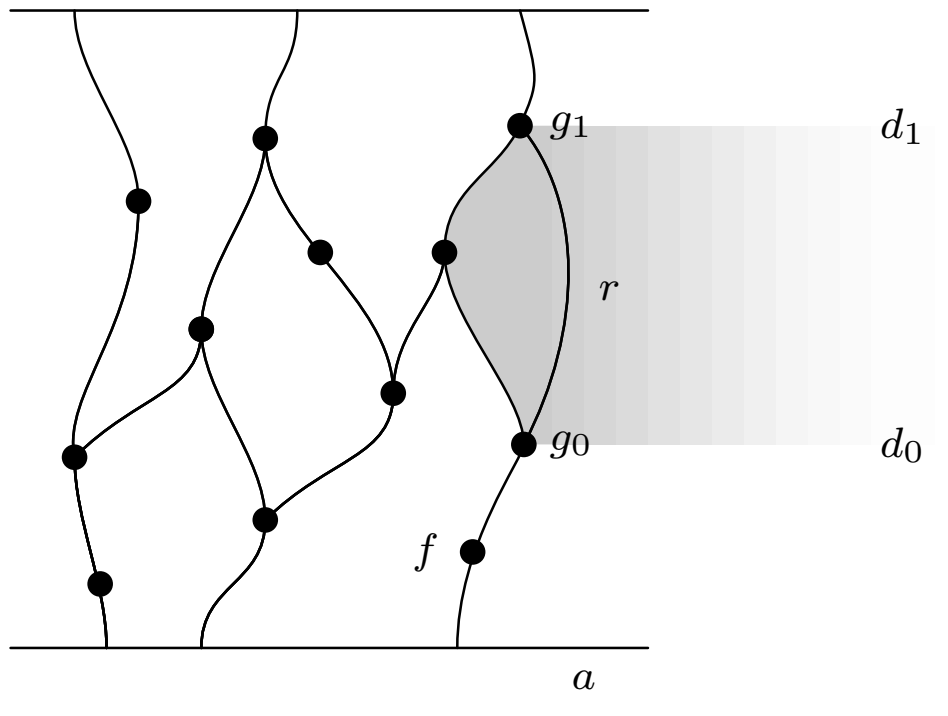

Now for any shape, start by erasing all $(1,1)$-valent nodes in the rightmost complete track. We claim that then a removable right edge exists. Remove this edge. Now we have reduced the number of edges, so by induction we arrive at the empty diagram.

To prove the claim, suppose the rightmost complete track contains a node $N$ of valence $(p, q)$ — otherwise it consists of a single edge which is clearly removable. Now either $p > 1$ or $q > 1$; without loss of generality we assume $q > 1$, so $N$ has more than one output edge. Now follow the rightmost output edge of $N$. If that edge $E$ ends at level $b_1$ or if the next node has more than one input edge, then $E$ is a removable right edge. Otherwise $E$ ends in a node of valence $(1, q')$ with $q' > 1$; then we can repeat the argument — since the graph is finite, eventually we come to a removable edge. □

Remark 2.5. The crucial condition for having an acyclic space of diagrams is that there are no 'floating islands' in the diagrams. Weaker conditions than being effective can preclude this, for example requiring only that every node has at least one output edge (yielding diagrams all of whose connected components are attached to the output line). The proof can easily be modified to cover such cases, but we do not need this.

**2.6. Free train track diagram categories.** If a tensor scheme $\mathscr{D}$ has the property that all its arrows have non-empty words as source and target, then it generates a free semimonoidal category, which is the category of isotopy classes of effective plane diagrams in $\mathscr{D}$. We shall use the following special case. Given a set $\mathscr{O}$, consider the tensor scheme $\mathscr{D}$ with $\mathscr{D}_0 = \{I\}$, and

$$\mathscr{D}_1 = \{f_o : I \to I \mid o \in \mathscr{O}\} \cup \{\alpha : II \to I\} \cup \{\beta : I \to II\}.$$

Since all the arrows have positive powers of $I$ as source and target, a free semimonoidal category is generated whose arrows are isotopy classes of effective diagrams in $\mathscr{D}$. Since $\alpha$ is of valence $(2, 1)$ and $\beta$ of valence $(1, 2)$, and all other generating arrows are of valence $(1, 1)$, the effective plane diagrams in $\mathscr{D}$ are precisely the train track diagrams in $\mathscr{D}$, the sort of diagrams drawn in Section 1.

## 3. From braided monoidal category to weak unit

**3.1. Outline of the idea.** Given any braided monoidal category $\mathbb{B}$, we are going to construct a semimonoidal 2-category $\mathscr{C}$ with weak unit $I$, such that $\mathrm{End}_{\mathscr{C}}(I)$ is equivalent to $\mathbb{B}$ as a braided monoidal category. The strategy is first to take the underlying semimonoidal category of $\mathscr{C}$ to be a free train track category, and construct a surjective set map from each hom set to the object set of $\mathbb{B}$. Then define the 2-cells in $\mathscr{C}$ by pulling back the 1-cells from $\mathbb{B}$. Hence by construction each hom cat of $\mathscr{C}$ is equivalent to $\mathbb{B}$. The braiding in $\mathrm{End}_{\mathscr{C}}(I)$ will correspond to the braiding in $\mathbb{B}$.

More specifically, let $\mathscr{O}$ be the object set of $\mathbb{B}$, and let $\mathbb{F}_{\mathrm{br}}(\mathscr{O})$ denote the free braided monoidal category on $\mathscr{O}$, with its projection to $\mathbb{B}$. The braided monoidal category $\mathbb{F}_{\mathrm{br}}(\mathscr{O})$ is naturally equivalent to the fundamental groupoid of the space of $\mathscr{O}$-labelled configurations of points in $\mathbb{R}^2$, and intuitively, from each hom set of the train track category there is a map to the configuration space, consisting in forgetting the tracks and only retain the trains. However, this map is not really well-defined, because the train track diagram category concerns only deformation classes of train track diagrams, not the actual diagrams. So to get a well-defined map we need to pass to the (equivalent) categories of *cliques* in $\mathbb{F}_{\mathrm{br}}(\mathscr{O})$ and $\mathbb{B}$.

In order to streamline the actual construction we shall first gather some facts about cliques and about configuration spaces and free braided monoidal categories.

## Preliminaries on cliques

**3.2. Cliques.** We shall recall some basic facts about cliques. (See Joyal-Street [**4**], p. 58.) Given a set (or a collection) $I$, let $\overline{I}$ denote the groupoid whose object set is $I$ and whose arrow set is $I \times I$ with the two projections as source and target. If $I$ is nonempty then $\overline{I}$ is contractible. A *clique* in a category $\mathscr{C}$ is a functor $\overline{I} \to \mathscr{C}$ for some nonempty set $I$. In other words, it consists of a nonempty collection of objects $(x_i \mid i \in I)$, together with a collection of comparison isomorphisms $x_{ij} : x_i \to x_j$ satisfying $x_{ii} = \mathrm{id}_{x_i}$ and $x_{ij} x_{jk} = x_{ik}$, for all $i, j, k \in I$. A *morphism* from clique $(x_i \mid i \in I)$ to clique $(y_j \mid j \in J)$ is a natural transformation

$$\begin{array}{ccc} \overline{I \times J} & \longrightarrow & \overline{J} \\ \downarrow & \nearrow & \downarrow{\scriptstyle y} \\ \overline{I} & \underset{x}{\longrightarrow} & \mathscr{C} \end{array}$$

In other words, it consists in a collection of arrows in $\mathscr{C}$

$$(f_{ij} : x_i \to y_j \mid i \in I, j \in J)$$

such that this square commutes:

$$\begin{array}{ccc} x_i & \xrightarrow{f_{ij}} & y_j \\ {\scriptstyle x_{ip}}\downarrow & & \downarrow{\scriptstyle y_{jq}} \\ x_p & \underset{f_{pq}}{\longrightarrow} & y_q \end{array}$$

Note that a morphism is completely determined by specifying any one of its components $f_{ij}$.

Let $\widetilde{\mathscr{C}}$ denote the category of cliques in $\mathscr{C}$. There is a canonical equivalence of categories $\mathscr{C} \to \widetilde{\mathscr{C}}$ given by sending an object $x \in \mathscr{C}$ to the singleton clique supported at $x$. (There is no canonical functor in the other direction.)

**3.3. Cliques in monoidal categories.** If $(\mathscr{C}, \otimes, \Bbbk)$ is a monoidal category, then there is a canonical monoidal structure on $\widetilde{\mathscr{C}}$: the tensor product $\widetilde{\otimes}$ is defined point-wise:

$$(x \widetilde{\otimes} y)_{i,j} := x_i \otimes y_j \qquad \text{and} \qquad (x \widetilde{\otimes} y)_{(i,j)(p,q)} := x_{ip} \otimes y_{jq},$$

the indexing set of the tensor product being $I \times J$. The neutral object $\widetilde{\Bbbk}$ is the singleton clique $* \mapsto \Bbbk$. (Note that even if $(\mathscr{C}, \otimes, \Bbbk)$ is a strict monoidal category, $(\widetilde{\mathscr{C}}, \widetilde{\otimes}, \widetilde{\Bbbk})$ will not be strict, since it involves the non-strictness of the cartesian products of the indexing sets.) If $(\mathscr{C}, \otimes, \Bbbk)$ has a braiding $\tau_{p,q} : p \otimes q \to q \otimes p$, then there is induced a braiding on $(\widetilde{\mathscr{C}}, \widetilde{\otimes}, \widetilde{\Bbbk})$ too: the components of $x \widetilde{\otimes} y \to y \widetilde{\otimes} x$ are simply $\tau_{x_i, y_j} : x_i \otimes y_j \to y_j \otimes x_i$.

**3.4. Lowershriek of a clique.** If $x : \overline{I} \to \mathscr{C}$ is a clique, and $F : \mathscr{C} \to \mathscr{D}$ is any functor, then obviously the composite $\overline{I} \to \mathscr{C} \to \mathscr{D}$ is again a clique in $\mathscr{D}$ which we denote by $F_!x$. This construction defines a functor $F_! : \widetilde{\mathscr{C}} \to \widetilde{\mathscr{D}}$. If $F : \mathscr{C} \to \mathscr{D}$ is a monoidal functor then there is induced a monoidal structure on $F_! : \widetilde{\mathscr{C}} \to \widetilde{\mathscr{D}}$ as well, and if $F$ is braided monoidal then so is $F_!$.

**3.5. Inverse image clique.** If $x : \overline{I} \to \mathscr{D}$ is a clique, and $F : \mathscr{C} \xrightarrow{\sim} \mathscr{D}$ is an equivalence of categories, then the 2-fibre product

$$F^*x := \mathscr{C} \underset{\mathscr{D}}{\times} \overline{I} \ \to \mathscr{C}$$

is a clique in $\mathscr{C}$. Specifically, for the 2-fibre product $F^*x$ we take the category whose objects are triples $(c, i, \gamma)$ where $c \in \mathscr{C}$, $i \in I$, and $\gamma : F(c) \xrightarrow{\sim} x_i$ is a specified comparison isomorphism, and whose arrows from $(c_0, i_0, \gamma_0)$ to $(c_1, i_1, \gamma_1)$ are arrows $\omega : c_0 \to c_1$ making this diagram commute:

$$\begin{array}{ccc} F(c_0) & \xrightarrow{\gamma_0} & x_{i_0} \\ {\scriptstyle F(\omega)}\downarrow & & \downarrow{\scriptstyle x_{i_0 i_1}} \\ F(c_1) & \xrightarrow[\gamma_1]{} & x_{i_1}. \end{array}$$

(Essential surjectivity of $F$ ensures that $F^*x$ is nonempty, and fully faithfulness ensures that $F^*x$ is contractible.) This construction defines a functor $F^* : \widetilde{\mathscr{D}} \to \widetilde{\mathscr{C}}$, which is again an equivalence of categories. If $F : \mathscr{C} \to \mathscr{D}$ is a monoidal functor then there is induced a monoidal structure on $F^* : \widetilde{\mathscr{D}} \to \widetilde{\mathscr{C}}$, and if $F$ is braided monoidal then so is $F^*$. (Note that even if $F$ is strictly monoidal, $F^*$ will not be strict except if $\mathscr{D}$ is discrete.)

## Configuration spaces and free braided monoidal categories

**3.6. Labelled configuration spaces.** Given a set $\mathscr{O}$, let $C_n(\mathbb{R}^2, \mathscr{O})$ denote the space of configurations of $n$ distinct points in $\mathbb{R}^2$, each labelled by an element in $\mathscr{O}$. We are interested in the disjoint union

$$C(\mathbb{R}^2, \mathscr{O}) := \coprod_{n \geq 0} C_n(\mathbb{R}^2, \mathscr{O}).$$

In other words, $C(\mathbb{R}^2, \mathscr{O})$ is the space of functions $S \to \mathscr{O}$ where $S$ is a finite subset of $\mathbb{R}^2$. When $\mathscr{O}$ is the singleton set, $C(\mathbb{R}^2, \mathscr{O})$ is the standard space of configurations of points in $\mathbb{R}^2$, and its fundamental groupoid $\Pi_1(C(\mathbb{R}^2))$ is equivalent to the braid category (the free braided monoidal category on one object).

If $\Lambda$ is any set, and $\rho : \Lambda \to C(\mathbb{R}^2, \mathscr{O})$ is a set map, we denote by $\Pi_1(C(\mathbb{R}^2, \mathscr{O}), \rho)$ the groupoid whose 0-cells are the elements of $\Lambda$, and whose 1-cells are pulled back from $\Pi_1 C(\mathbb{R}^2, \mathscr{O})$. That is, if $x$ and $y$ are elements in $\Lambda$, then $\operatorname{Hom}(x, y) := \operatorname{Hom}_{\Pi_1 C(\mathbb{R}^2, \mathscr{O})}(\rho(x), \rho(y))$. By construction there is a fully faithful functor $\Pi_1(C(\mathbb{R}^2, \mathscr{O}), \rho) \to \Pi_1 C(\mathbb{R}^2, \mathscr{O})$, which is an equivalence of categories provided $\rho$ is surjective on the set of connected components of $\Pi_1 C(\mathbb{R}^2, \mathscr{O})$.

Given a set $\mathscr{O}$, let $\mathbb{F}_{\mathrm{br}}(\mathscr{O})$ denote the free braided monoidal category on $\mathscr{O}$: the underlying monoid of $\mathbb{F}_{\mathrm{br}}(\mathscr{O})$ is $\mathscr{O}^*$, the free monoid on $\mathscr{O}$, and the arrows are the $\mathscr{O}$-coloured braids. For each word $\mathbf{v} = v_1 \cdots v_n$ in $\mathscr{O}^*$ consider the labelled configuration supported on the positive-integer points of the $x$-axis, where the point $(i, 0)$ has label $v_i$. This defines a set map $\rho : \mathscr{O}^* \to C(\mathbb{R}^2, \mathscr{O})$, and the free braided monoidal category $\mathbb{F}_{\mathrm{br}}(\mathscr{O})$ is naturally isomorphic to $\Pi_1(C(\mathbb{R}^2, \mathscr{O}), \rho)$.

Using integer points is in fact an arbitrary choice, and intuitively we are really talking about cliques: for each word $\mathbf{v} \in \mathscr{O}^*$, there is a clique whose objects are ordered $\mathbf{v}$-labelled configurations on the $x$-axis $\mathbb{R}^1$, and whose connecting isomorphisms are the (homotopy classes of) order-preserving paths in $C(\mathbb{R}^1, \mathscr{O})$ — the positive-integer point configurations are just normal-form representatives for these cliques. When making drawings, the extra flexibility is important.

Combining these two functors we get an equivalence of categories

$$\mathbb{F}_{\mathrm{br}}(\mathscr{O}) \xrightarrow{\sim} \Pi_1 C(\mathbb{R}^2, \mathscr{O}).$$

This is the geometric model of $\mathbb{F}_{\mathrm{br}}(\mathscr{O})$. The tensor operation is just concatenation of configurations on the $x$-axis; the braiding is the homotopy class of the movement whereby some points on the right move up in the upper halfplane and left past some points on the left (and back to the $x$-axis).

## The construction

**3.7. The 1-skeleton of $\mathscr{C}$ as free semi-monoidal category.** Let $(\mathbb{B}, \otimes, \Bbbk, \tau)$ be a braided monoidal category with object set $\mathscr{O}$. We are going to construct a strict semimonoidal 2-category with weak unit. Since this category is going to have strict composition laws and strict tensor product, it makes sense first to construct its 1-skeleton, a semimonoidal category, and then describe the 2-cells afterwards.

Let $\mathscr{C}_0$ denote the free train track category on $\mathscr{O}$, as in 2.6. The object set is $\{I, I^2, I^3, \dots\}$ and the arrows are isotopy classes of train track diagrams with trains labelled in $\mathscr{O}$.

We employ the notation $[I^p, I^q]_0$ for the hom sets in this category, and we are going to enrich over ***Cat*** to arrive at the hom cats $[I^p, I^q]$ defining the 2-category $\mathscr{C}$. The object $I$ is going to be a weak unit for $\mathscr{C}$, but note that it is not a unit in $\mathscr{C}_0$, because $\alpha : II \to I$ is not an isomorphism.

Let $X_{\mathscr{O}}$ denote the space of all train track diagrams labelled in $\mathscr{O}$, and let $X_{\mathscr{O}}(p, q)$ denote the subspace of train track diagrams with $p$ inputs and $q$ outputs. The set $[I^p, I^q]_0$ is just the set of connected components of $X_{\mathscr{O}}(p, q)$.

Let $\epsilon : X_{\mathscr{O}} \to C(\mathbb{R}^2, \mathscr{O})$ denote the map that forgets the tracks and only retains the trains, i.e., returns the configuration of the $\mathscr{O}$-labelled nodes of a given diagram. (Note that this maps also forgets the positions of $\alpha$- and $\beta$-labelled nodes.)

Passing to the fundamental groupoid of these spaces we have the following diagram:

$$\begin{array}{ccccc}
\Pi_1 X_{\mathscr{O}}(p,q) & \xrightarrow{\ \epsilon\ } & \Pi_1 C(\mathbb{R}^2, \mathscr{O}) & & \\
{\scriptstyle \pi_0}\downarrow & & \uparrow{\scriptstyle \rho} & & \\
[I^p, I^q]_0 & & \mathbb{F}_{\mathrm{br}}(\mathscr{O}) & \xrightarrow{\ \gamma\ } & \mathbb{B}
\end{array}$$

Since $\pi_0$ has contractible fibres by Proposition 2.4, each element in the set $[I^p, I^q]_0$ defines a clique in $\Pi_1 X_{\mathscr{O}}(p, q)$ namely the inclusion of the fibre. Composing with $\epsilon_!$ and $\rho^*$, and finally with $\gamma_!$, we get a map denoted $\theta$:

$$[I^p, I^q]_0 \longrightarrow \widetilde{\Pi_1 X}_{\mathscr{O}}(p,q) \xrightarrow{\epsilon_!} \widetilde{\Pi_1 C}(\mathbb{R}^2, \mathscr{O}) \xrightarrow{\rho^*} \widetilde{\mathbb{F}_{\mathrm{br}}(\mathscr{O})} \xrightarrow{\gamma_!} \widetilde{\mathbb{B}}$$

Given an element $D \in [I^p, I^q]_0$, i.e. an isotopy class of train track diagrams, let us trace through the clique maps to get a more concrete description of the clique $\theta D$. In the following, the word 'generically' means that the involved configurations are assumed not to intersect — this assumption is convenient for the sake of drawing pictures.

The image clique in $\Pi_1 C(\mathbb{R}^2, \mathscr{O})$ has as objects those labelled configurations that can support a train track diagram of isotopy class $D$. Note that the progressive condition on the diagrams implies a restriction on the possible configurations: if $a$ is a dot in $D$ that comes before a dot $b$ on the same complete track, then clearly this order must be reflected in the $y$-coordinate of the corresponding points in the configuration. The comparison arrows in the clique are homotopy classes of paths in $C(\mathbb{R}^2, \mathscr{O})$ such that every intermediate configuration can also support a train diagram of class $D$. If one concrete configuration $F_0$ (being the train points of some diagram $D_0$) is chosen as representative for the clique, then generically another representing configuration $F_1$ together with the connecting isomorphism $\phi : F_0 \to F_1$) can be depicted as another configuration connected to $F_0$ with non-intersecting parallel strings.

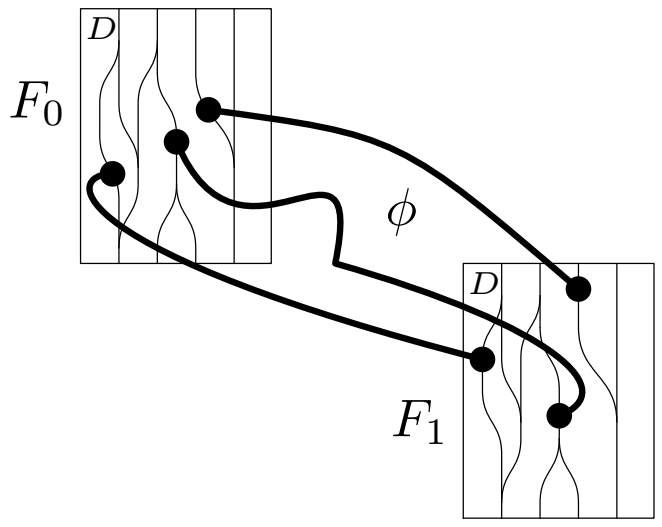

The rigidity of train track diagrams expressed by Proposition 2.4 means that we essentially can regard each configuration as a black box, and the connecting isomorphisms are essentially just translations of such boxes.

The lowershriek of this clique in $\widetilde{\mathbb{F}_{\mathrm{br}}(\mathscr{O})}$ has this description: the objects are triples $(W, F, \gamma)$ where $W$ is a labelled configuration on the $x$-axis, $F$ is a configuration that can support a train track diagram of class $D$, and $\gamma$ is a homotopy class of paths from $W$ to $F$, which we can think of as a linearisation of the set of train points. We depict the triple $(W, F, \gamma)$ as a *string configuration* like this:

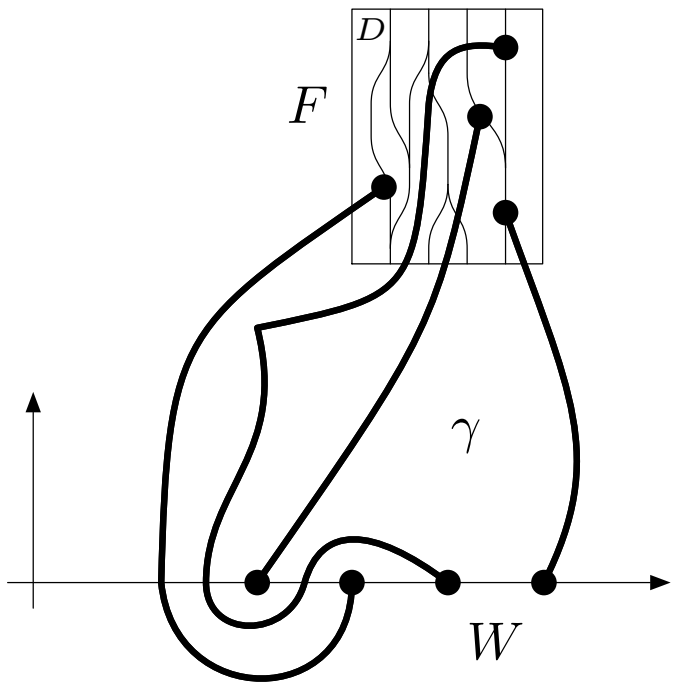

The connecting isomorphism from $(W_0, F_0, \gamma_0)$ to $(W_1, F_1, \gamma_1)$ is the homotopy class of paths from $W_0$ to $W_1$ compatible with $\gamma_0$, $\gamma_1$, and $\phi$ (where $\phi : F_0 \to F_1$ is the connecting isomorphism in the clique of configurations corresponding to $D$, as described above). For simplicity we assume that $F_0 = F_1$, then generically the connecting isomorphism $\omega$ is obtained by drawing non-intersecting strings from $W_0$ to $W_1$ in the complement of the strings representing $\gamma_0$ and $\gamma_1$.

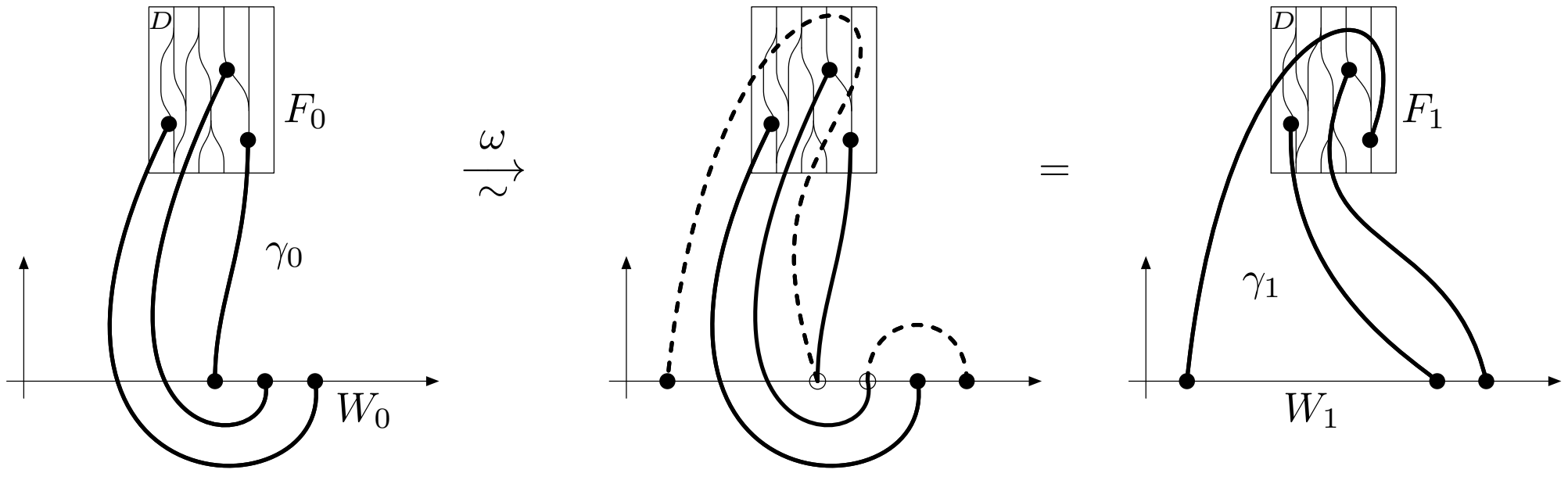

It is convenient to choose representing string configurations in such a way that $F$ as well as all the strings are contained in the upper half-plane.

**3.8. The composition law in $\mathscr{C}$ on the object level.** $\mathscr{C}_0$ is a semimonoidal category by construction. We shall briefly indicate the interpretation of its structures on the level of cliques. The tensor product operation on train track diagrams,

$$[I^p, I^q]_0 \times [I^m, I^n]_0 \xrightarrow{\otimes} [I^{p+m}, I^{q+n}]_0$$

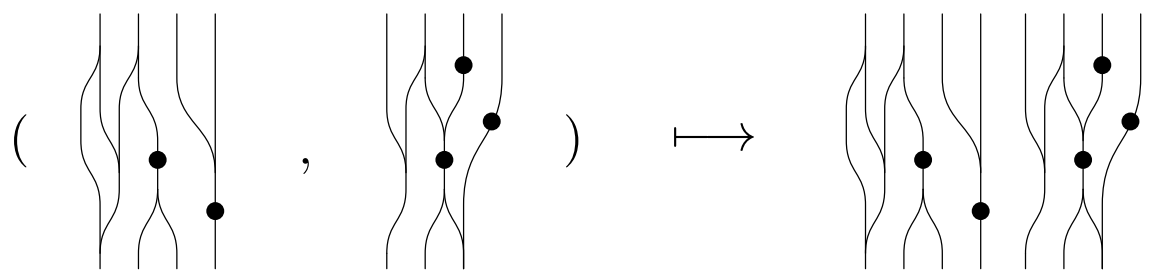

is just horizontal concatenation of cliques of string configurations, and hence corresponds to the tensor product in $\widetilde{\mathbb{B}}$.

The composition law

$$[I^p, I^q]_0 \times [I^q, I^r]_0 \longrightarrow [I^p, I^r]_0$$

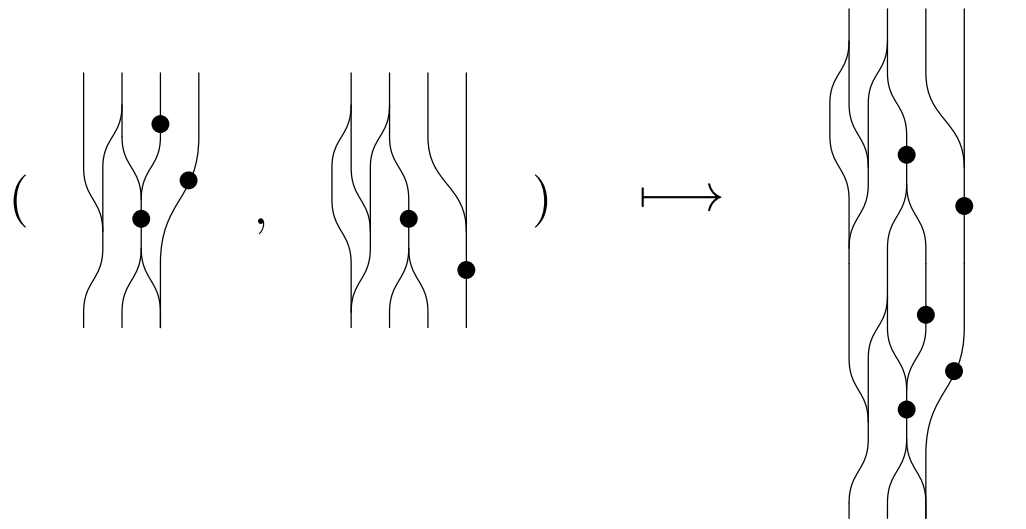

corresponds in $\widetilde{C(\mathbb{R}^2, \mathscr{O})}$ to vertical stacking of cliques of configurations of points. In $\widetilde{\mathbb{F}_{\mathrm{br}}(\mathscr{O})}$, the picture for this operation is this

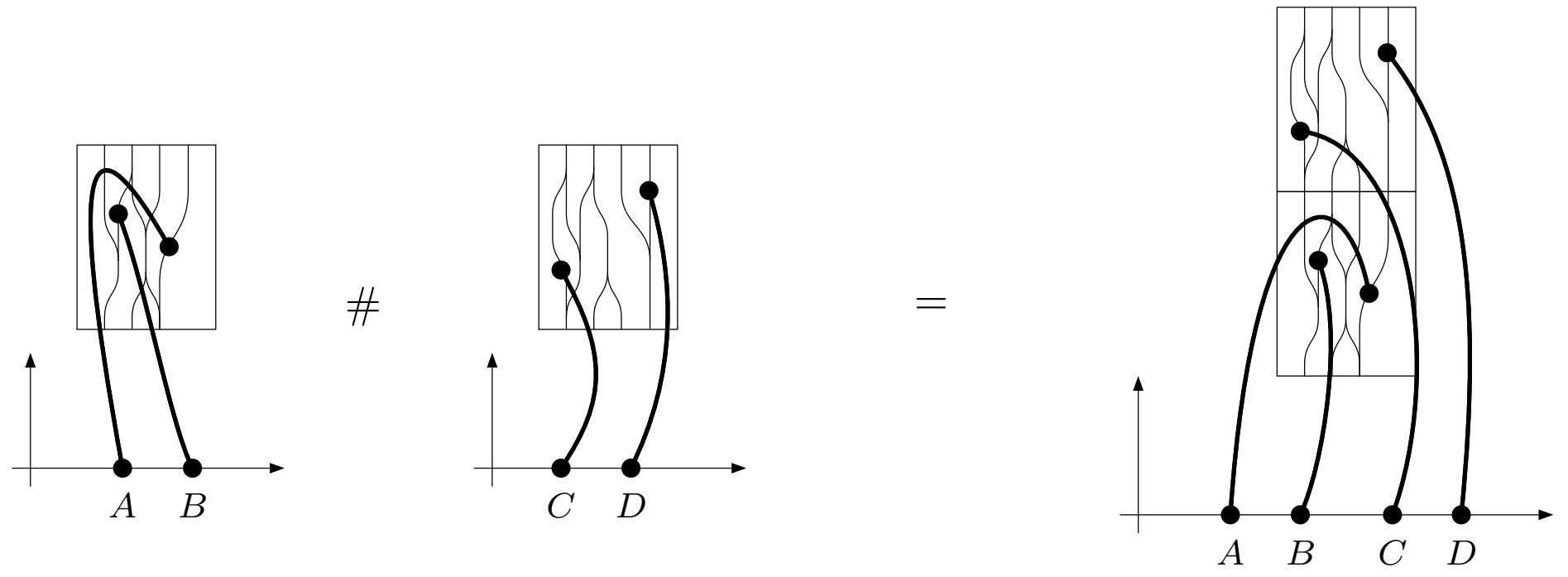

Hence it corresponds to concatenation of words in $\widetilde{\mathbb{F}_{\mathrm{br}}(\mathscr{O})}$, and down in $\widetilde{\mathbb{B}}$ it corresponds to the tensor product. However this is only one of the possible representatives for the clique corresponding to the composite; the opposite ordering of the factors constitutes another representative, and the unique connecting isomorphism corresponds to the braiding, as indicated here:

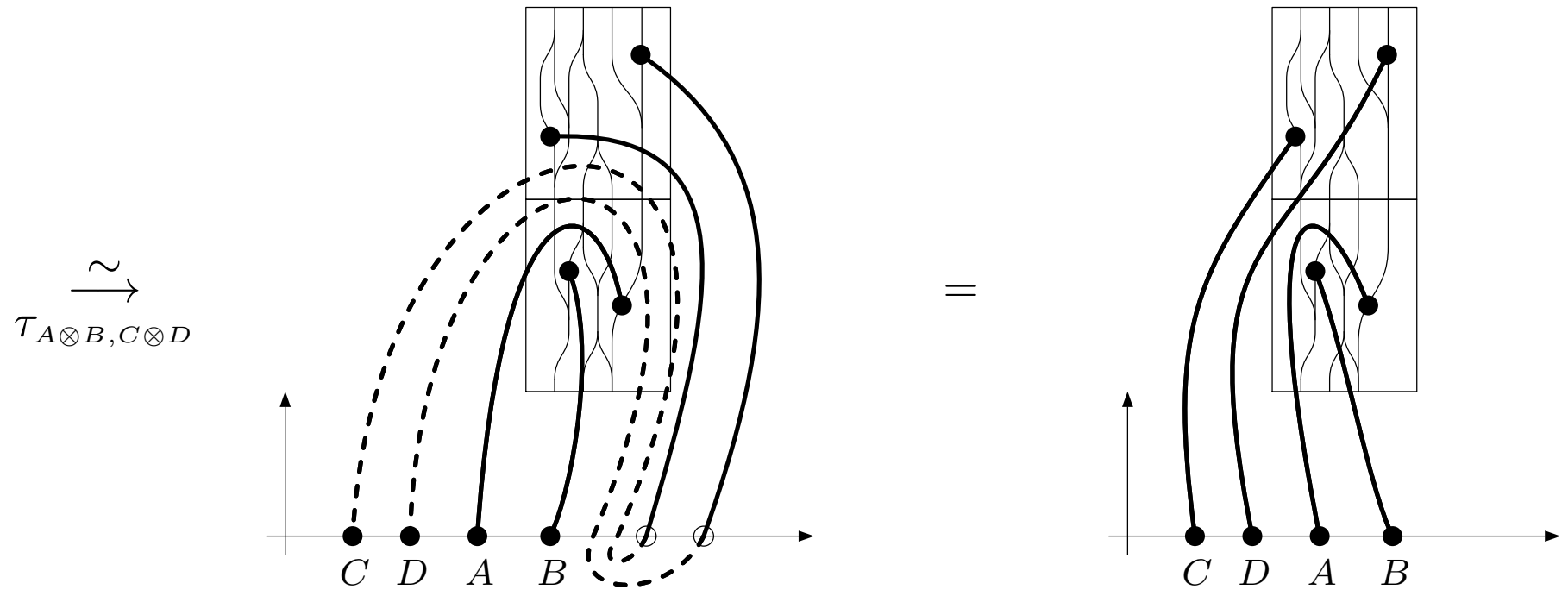

Finally we explain the interchange law in terms of cliques of string configurations. Given four train track diagrams classes, composable as indicated:

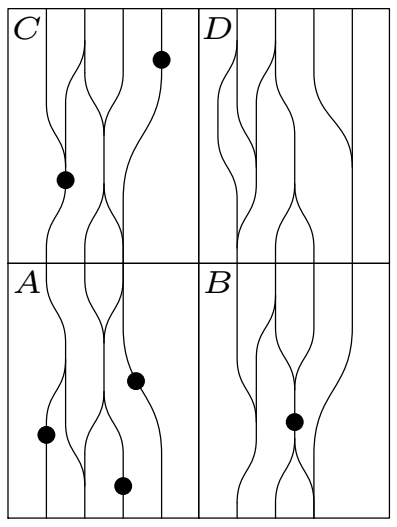

the interchange law reads

$$(A \otimes B) \# (C \otimes D) = (A \# C) \otimes (B \# D).$$

Taking representatives for the composites as above, this equation has the following interpretation in $\widetilde{\mathbb{F}_{\mathrm{br}}(\mathscr{O})}$ (or in $\widetilde{\mathbb{B}}$):

(3)

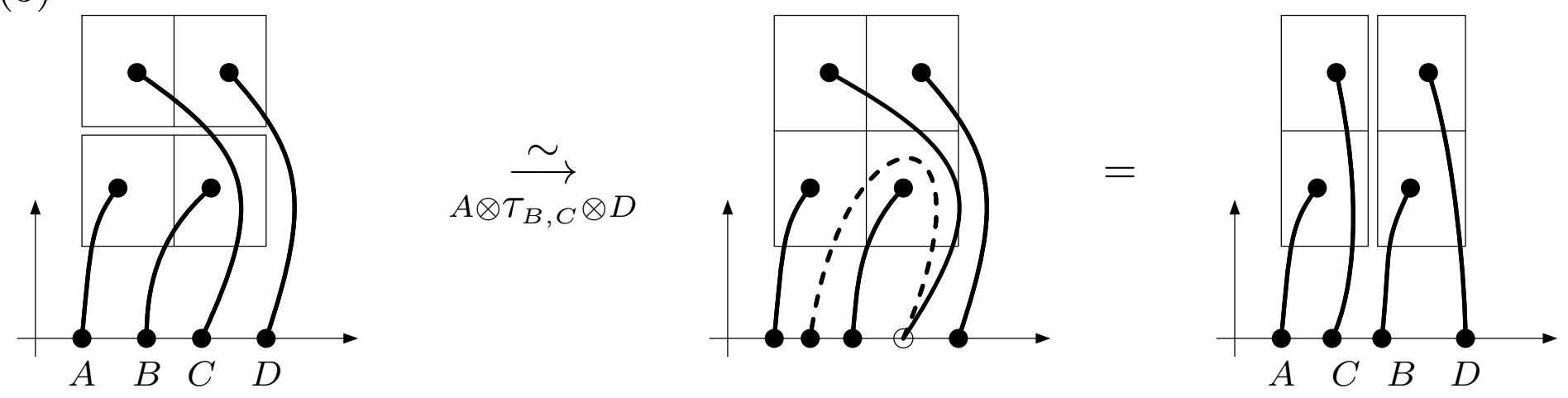

Again, the unique isomorphism connecting the representatives for the total composite is just an instance of the braiding in $\mathbb{B}$.

**3.9. The 2-cells of $\mathscr{C}$.** The 2-cells of $\mathscr{C}$ are defined by pulling back the 1-cells in $\widetilde{\mathbb{B}}$ along $\theta$:

$$[D, D'] := \widetilde{\mathbb{B}}(\theta D, \theta D').$$

A clique map $\theta D \to \theta E$ is represented by a 1-cell $f$ in $\mathbb{B}$ between the chosen representatives for $\theta D$ and $\theta E$. In the following drawing the two planes picture representatives for the $\theta D$ and $\theta E$, and $f$ is indicated as the grey graph between the bottom lines. (Although this graph looks planar, it is meant as a 3D diagram like in [**4**], Chapter 3.)

(4)

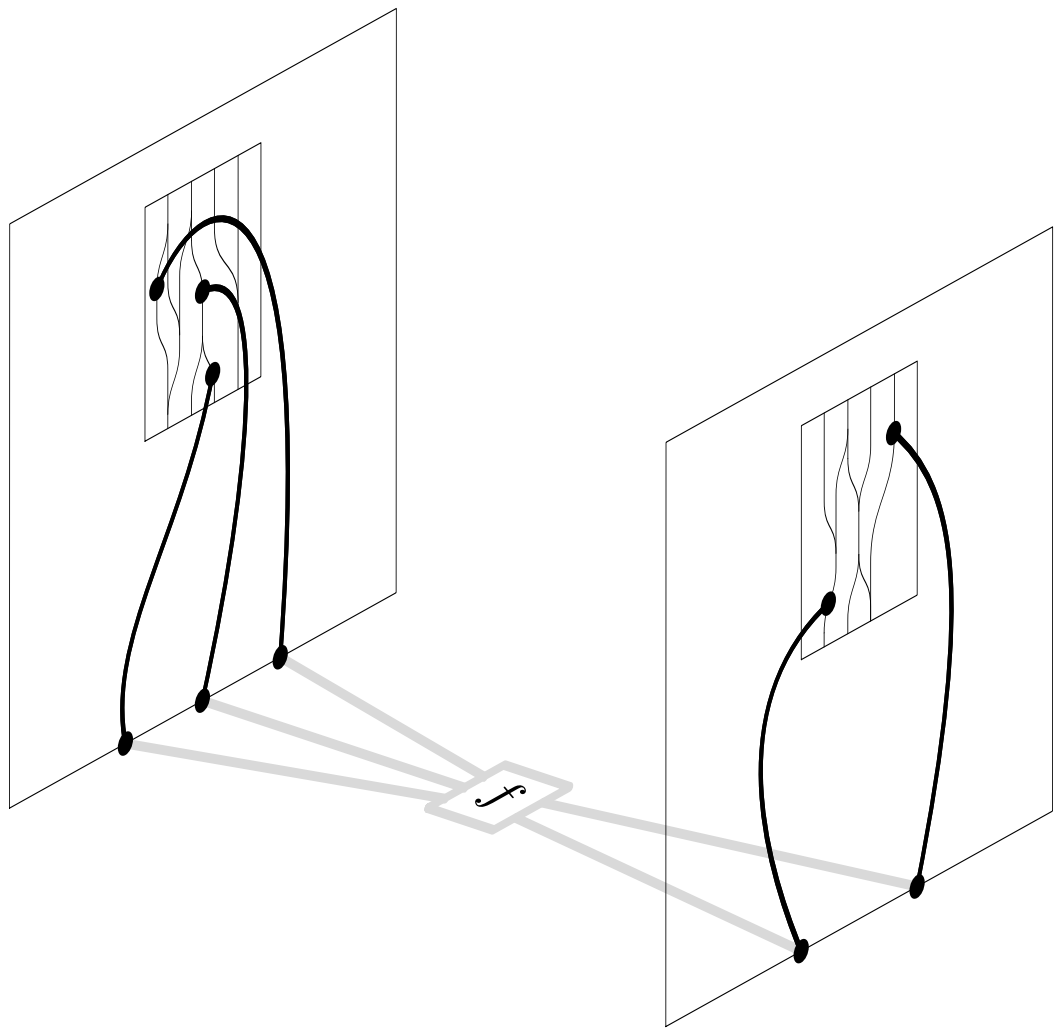

The different representatives are uniquely related by conjugation with components of the braiding $\tau$. These are given as part of the clique structure, but they can also be characterised in terms of the figures as those braidings that can be realised in the complement of the strings representing the linearisations. In other words, if $A$ and $A'$ are two representatives for the clique $\theta D$, connected by $u : A \to A'$, and if $v : B \to B'$ is a connected pair of representatives for $\theta E$, then $\mathbb{B}$-arrows $f : A \to B$ and $f' : A' \to B'$ represent the same clique map when this square commutes:

$$\begin{array}{ccc} A & \xrightarrow{f} & B \\ {\scriptstyle u}\downarrow & & \downarrow{\scriptstyle v} \\ A' & \xrightarrow[f']{} & B' \end{array}$$

**3.10. Composition law and tensor product.** The set maps $[I^p, I^q]_0 \times [I^q, I^r]_0 \to [I^p, I^r]_0$ extend to functors by defining the horizontal composition of 2-cells in $\mathscr{C}$ to be the tensor product of the representing 1-cells down in $\widetilde{\mathbb{B}}$. Different ways of writing the representing tensor product are uniquely related by isomorphisms, and these are just components of the braiding $\tau$. The 'vertical' composition of 2-cells is just composition of arrows in $\widetilde{\mathbb{B}}$.

Finally, the tensor product on $\mathscr{C}$, paralleling train track diagrams, extends in the same way to the new 2-cells. Functoriality, i.e., the interchange law on the level of 2-cells, follows from the same argument as in Figure (3), just applying the braiding $\tau$ to arrows instead of merely to objects.

This concludes the description of the semimonoidal 2-category $\mathscr{C}$.

**3.11. Weak unit.** Now we have constructed a semimonoidal 2-category with an object $I$ such that $\mathrm{End}(I)$ is equivalent to $\mathbb{B}$. We now check that $I$ is a weak unit. First, the arrow $\alpha : II \to I$ is an equi-arrow in $\mathscr{C}$, with quasi-inverse $\beta : I \to II$. Indeed, the required invertible 2-cells can both be represented in $\widetilde{\mathbb{B}}$ by the identity arrow of the unit object $\Bbbk$ — in terms of pictures like Fig 4 this is just the empty braid between empty configurations on the bottom lines.

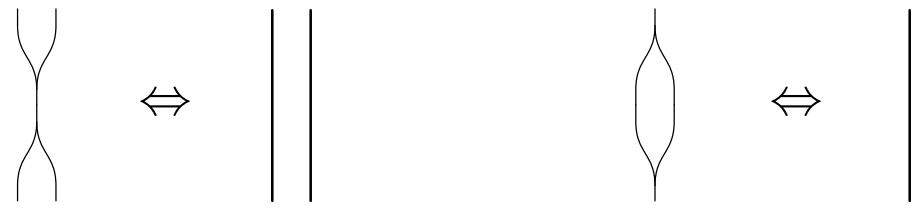

We then have to check that tensoring with $I$ on the right (or on the left) is an equivalence of 2-categories $\mathscr{C} \to \mathscr{C}$. In other words,

$$[I^p, I^q] \quad \longrightarrow \quad [I^{p+1}, I^{q+1}]$$

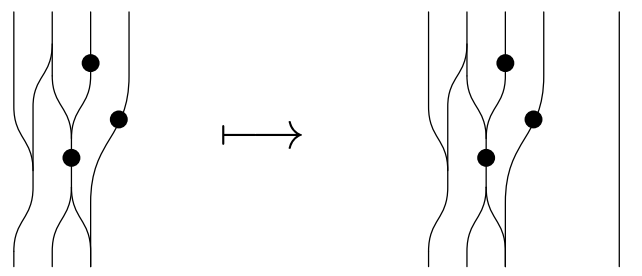

should be an equivalence of categories. This functor is fully faithful: given two diagrams

$$I^p \underset{D'}{\overset{D}{\rightrightarrows}} I^q$$

the 2-cells $D \Rightarrow D'$ are given by $\mathbb{B}$-arrows placed between the bottom lines of two representing diagrams like in figure 4. The description is exactly the same for the 2-cells $DI \Rightarrow D'I$ because the extra trainless track on the right doesn't show up on the bottom lines. Finally, the functor is essentially surjective: any diagram $E \in [I^{p+1}, I^{q+1}]$ is isomorphic to one of the form $DI$ for $D \in [I^p, I^q]$ — take any train track diagram with the same set of trains as $E$, then the corresponding cliques are both represented by tensor products in $\mathbb{B}$ with the same factors, so a suitable braid between the two tensor products provides the desired isomorphism.

It is clear that the braiding on $\mathrm{End}(I)$ constructed in Section 1 corresponds to the braiding in $\mathbb{B}$.

This concludes the proof of the Main Theorem. □

Remark 3.12. In a sense, the key point of the construction is that the two non-strict tensor products on $C(\mathbb{R}^2, \mathscr{O})$ are strictified by quotienting by an equivalence relation, consisting in not caring about the precise position of the points but only their relative position. In order to keep track of this relative position, some grid or background texture is needed, to prevent the points from moving around each other (which would lead to the Eckmann-Hilton argument), and introducing this grid gives rise to the weak units (which are pure grid, no points). This background grid itself must be sufficiently rigid and attached. This is achieved by excluding the unit object in the horizontal direction (the excluded $I^0$); Proposition 2.4 is a formal expression of this idea. For the same reason, it is necessary to give up the vertical monoidal structure for a many-object version expressed by the variable number of strings.

A different approach to such strictifying grids, based on subdivided rectangles instead of train tracks, was presented by the second named author at the conference on *Higher-Order Geometry and Categorification* in Lisbon, July 2003. Indeed that method does strictify the two monoidal structures without breaking the interchange law, but such grids are not sufficiently rigid to prevent the braiding from collapsing to a symmetry. It is shown in [**7**] that such collapse will always happen in the 2-monoidal case, hence the necessity to replace one of the monoidal structures by a many-object version (but still contractible).

**3.13. Braided categorical groups.** If $\mathbb{B}$ is a braided categorical group, i.e. a monoidal groupoid such that every object has a monoidal inverse (up to isomorphism), then the corresponding $\mathscr{C}$ as in the construction above will clearly be a semimonoidal strict 2-groupoid whose tensor product is invertible up to equivalence with respect to the weak unit $I$. Conversely, for any such semimonoidal 2-groupoid with weak unit $I$, the braided monoidal category $\mathrm{End}_{\mathscr{C}}(I)$ will in fact be a braided categorical group. Since braided categorical groups are models for connected, simply connected homotopy 3-types, we have shown the Main Corollary, stated in the introduction.

## References

[1] Robert Gordon, A. John Power, and Ross Street. *Coherence for tricategories.* Mem. Amer. Math. Soc. **117** (1995), vi+81.

[2] Alexander Grothendieck. Pursuing stacks. Letter to D. Quillen, 1983, 600pp.

[3] André Joyal and Joachim Kock. *Coherence for weak units.* Manuscript in preparation.

[4] André Joyal and Ross Street. *The geometry of tensor calculus. I.* Adv. Math. **88** (1991), 55–112.

[5] G. Max Kelly and Ross Street. *Review of the elements of* 2*-categories.* In *Category Seminar (Proc. Sem., Sydney, 1972/1973)*, pp. 75–103. Lecture Notes in Math., Vol. 420. Springer, Berlin, 1974.

[6] Joachim Kock. *Weak identity arrows in higher categories.* Internat. Math. Res. Papers, vol. 2006, 1–54, (math.CT/0507116).

[7] Joachim Kock. *Note on commutativity in double semigroups and two-fold monoidal categories.* Preprint, math.CT/0608452.

[8] Carlos Simpson. *Homotopy types of strict 3-groupoids.* Preprint, math.CT/9810059.
[9] Zouhair Tamsamani. *Sur des notions de n-catégorie et n-groupoïde non strictes via des ensembles multi-simpliciaux.* $K$-Theory **16** (1999), 51–99. (alg-geom/9512006 and alg-geom/9607010.)

Département de mathématiques, Université du Québec à Montréal, Case postale 8888, succursale centre-ville Montréal (Québec), H3C 3P8 Canada